**Lagrange's Theorem On The Minimal Set Of Squares**
**N. A. Carella, April 2011.**

*Abstract*: It will be demonstrated that there is a thin basis of order four of minimal cardinality $\#A(x) = O(x^{1/4})$. The current literature shows the existence of a thin basis of order four of cardinality $\#A(x) = O(x^{1/4+\varepsilon})$, $\varepsilon > 0$, but speculates on the existence of a thin basis of order four of minimal cardinality.

## 1. Integers As Sums Of Squares

Let $\mathbb{N} = \{ 0, 1, 2, 3, \ldots \}$ denotes the set of nonnegative integers. A subset of nonnegative integers $A \subset \mathbb{N}$ is called a basis of order $s \geq 1$ if every nonnegative integer $n \in \mathbb{N}$ can be written as a sum of $s$ elements of the subset $A$, repetitions are permitted. The subset of squares $A = \{ 0, 1, 4, 9, \ldots \}$ is the basic prototype of a basis of the integers. It has been known for quite sometime that the subset of squares forms a basis of order four.

Let $x \geq 1$ be a real number and let $A(x) = \{ n^2 \leq x : 0 \leq n \leq x^{1/2} + 1 \} \subset A$ be the subset of squares up to $x$. Since there are $x^{1/2}$ squares in the interval $[1, x]$, the maximal size of a basis of order four is $\#A(x) \leq x^{1/2} + 1$. On the other hand, if $A(x)$ is a basis of order four, then the total number of four-tuples $a_1 + a_2 + a_3 + a_4 \leq x$ with $a_i \in A(x)$ must satisfy $\#A(x)^4 \geq c_0 x$, $c_0 \geq 1$ constant. Accordingly, the minimal size of a basis of order four is $\#A(x) \geq c_1 x^{1/4}$, $c_1 \geq 1$ constant.

The existence of a thin basis of order four of cardinality $\#A(x) = O(x^{1/4+\varepsilon})$ of the set of integers was proved in [ES] and [ZR]. Further, it is unknown if the parameter $\varepsilon > 0$ can be removed, see [NT, p. 34] for a discussion.

In this note it will be shown that there is a thin basis of order four of minimal cardinality $\#A(x) = c_1 x^{1/4}$, $c_1 \geq 1$ constant. Moreover, the proof is significantly simpler than the previous related proofs of this problem and easily generalizes to other related Warring problems. The precise result is as follows.

***Theorem* 1.**  There exists a thin basis of squares of the set of nonnegative integers of order four and minimal cardinality $\#A(x) = cx^{1/4}$, $c \geq 1$ constant.

The proof appears in Section 3.

## 2. Basic Results In Sums Of Squares
The quintessential results on the theory of integer representations as sums of squares are the following.

***Theorem* 2.**  (Girard) Every integer $n = AB^2$ such that $A$ is squarefree and has no prime divisor $p = 4a + 3$ is the sum of two squares.

According to [DN], Girard stated this problem, and then Fermat restated it. But, its first proof by infinite descend, was given by Euler about a century later. Nevertheless, this result is usually credited to Fermat.

***Theorem 3.*** (Gauss) Every integer $n \ne 4^r(8s + 7)$ is the sum of three squares.

***Theorem 4.*** (Lagrange) Every integer is the sum of four squares.

Background materials and proofs are given in [GR, p. 273], [MW], et cetera. It is quite easy to show that Theorems 2 and 3 imply Theorem 4. Moreover, it is quite possible that Theorem 2 implies both Theorem 3 and Theorem 4.

The general concept considered here for estimating the cardinality of a thin basis of $k$th powers is captured in the following result.

***Lemma 5.*** Let $A(x)$ be a set of cardinality $\#A(x) = x^\alpha$, $\alpha > 0$, and let $sA(x) = \{\, a_1 + \cdots + a_s : a_i \in A(x) \,\}$ be the sum set of $s \ge 1$ copies of $A(x)$. Then the following hold.
(i) The number of sums $a_1 + a_2 + \cdots + a_s$ with repetitions and $a_i \in A(x)$ is $\#(sA(x)) \le x^{\alpha s}/s!$.
(ii) If $a_1 + a_2 + \cdots + a_s \le x$, then $\#A(x) \ge c x^{1/s}$, $c > 0$ constant.
(iii) If each $a_i \in A(x)$ is a $k$th power, then $\#A(x) \le c x^{1/k}$, $c > 0$ constant.

This Lemma shows that the cardinality $\#A(x)$ of a basis $A(x)$ of $k$th powers and order $s \ge 2$ is in the range $c_1 x^{1/s} \le \#A(x) \le c_2 x^{1/k}$, where $c_1, c_2 > 0$ are constants.

## 3. Proof Of The Theorem

For a real number $x \ge 1$, define the subsets of squares

$$A_1(x) = \{\, n^2 \le c_0 x^{1/2} : 0 \le n \le c_0^{1/2} x^{1/4} \,\}, \text{ and } A_2(x) = \{\, n^2 \le x : c_0^{1/2} x^{1/4} \le n \le c_1^{1/2} x^{1/2} \,\}, \qquad (1)$$

where $c_0, c_1 \ge 1$ are constants. The basic idea of the proof of Theorem 1 is to show that every integer $n \le x$ can be written as a sum of four squares as in

$$n = a_1 + a_2 + a_3 + a_4 \le x, \qquad (2)$$

where $a_1, a_2, a_3 \in A_1(x)$ and $a_4 \in A_2(x)$. This basic principle emanates from the fact that every integer $n \in \mathbb{N}$ has a representation of the form

$$n = [n^{1/2}]^2 + m, \qquad (3)$$

where $0 \le |m| \le 2\sqrt{n} + 1$, and the bracket $[\,x\,]$ denotes the largest integer function. The distribution of the fractional parts $((\sqrt{n})) = n - [\sqrt{n}]$ is a very interesting and difficult problem. It is known that $((\sqrt{n})) < 1/\sqrt{n}$ for infinitely many integers $n \in \mathbb{N}$, but it is unknown if $((\sqrt{n})) < 1/\sqrt[4]{n}$, see [HN] for advance material, but it is not required in this work. The simple representations (3) of the integers $n \ge 1$, and Theorem 3 reduce the complexity of the proof of the existence of a minimal thin basis to a few lines. Moreover, this idea easily generalizes to the other related cases of the Warring problem $n = x_1^k + x_2^k + \cdots + x_s^k$, where $k \ge 1$, and $s \ge 2$ are fixed parameters.

***Theorem* 1.** There exists a thin basis of squares of the set of nonnegative integers of order four and minimal cardinality $\#A(x) = cx^{1/4}$, $c \geq 1$ constant.

Proof: Let $n \in \mathbb{N}$, and consider the integer $m = n - [n^{1/2} - t]^2$, where $t = 0, 1, 2, \ldots$, is a small integer. It is immediate that

$$0 < n - [n^{1/2} - t]^2 = m \leq c_2 \sqrt{n}, \qquad (4)$$

where $c_2 \geq 2$ is a constant. Put $a_4 = [n^{1/2} - t]^2 \in A_2(x)$. Further, since the congruence equation $n - [n^{1/2} - t]^2 = m \equiv 0, 1, \pm 2, \pm 3, 4 \bmod 8$ holds for some $t = 0, 1$, the integer $m \geq 0$ has a representation as a sum of three squares

$$m = a_1 + a_2 + a_3 \leq c_2 x^{1/2}, \qquad (5)$$

where $a_1, a_2, a_3 \in A_1(x) = \{ n^2 \leq c_0 x^{1/2} : 0 \leq n \leq c_0^{1/2} x^{1/4} \}$, see Theorem 3. Solving for $n$ returns $n = [n^{1/2} - t]^2 + m = a_1 + a_2 + a_3 + a_4$, where $a_4 \in A_2(x) = \{ n^2 \leq x : c_0^{1/2} x^{1/4} \leq n \leq c_1^{1/2} x^{1/2} \}$. To complete the verification, notice that the cardinality of the subset of squares $A(x) = A_1(x) \cup A_2(x)$ is $\#A(x) \leq \#A_1(x) + \#A_2(x) = cx^{1/4}$. ∎

The next interesting case seems to be a thin basis of cubes of order $s = 7$ for the set of large integers $\{ n_0 \leq n = x_1^3 + x_1^3 + \cdots + x_7^3 : x_i \geq 0 \}$, where $n_0$ is fixed. It is apparent that the same technique, a greedy-type algorithm as above, can be utilized to confirm the existence of a thin basis of cubes of the set of nonnegative integers of order seven and minimal cardinality $\#A(x) = cx^{1/7}$, $c \geq 1$ constant. Lemma 5 predicts the existence of a thin basis of cubes $B(x) = \{ n^3 \leq x : 0 \leq n \leq cx^{1/3} \}$ of cardinality in the range $c_3 x^{1/7} \leq \#B(x) \leq c_4 x^{1/3}$

## 4. Primes as Sums of Squares
Theorems 1, 2 and 3 provide a complete classification of the prime numbers as sums of squares.

***Corollary* 6.** For every prime $p$ there is a representation of minimal number of squares up to signs and permutations of the fourtuple $x, y, z$, and $w$. The four residues classes are:
(i) A prime $p = 8n + 1 = x^2 + y^2$, where $x \neq 0$, and $y \neq 0$.
(ii) A prime $p = 8n + 3 = x^2 + y^2 + z^2$, where $x \neq 0$, $y \neq 0$, and $z \neq 0$.
(iii) A prime $p = 8n + 5 = x^2 + y^2$, where $x \neq 0$, and $y \neq 0$.
(iv) A prime $p = 8n + 7 = x^2 + y^2 + z^2 + w^2$, where $x \neq 0$, $y \neq 0$, $z \neq 0$, and $w \neq 0$.

The proofs of (i) and (iii), which are widely available in the literature, follow from Theorem 2. The proof of (ii) follows from Theorem 3. The proof of (iv) follow from Theorem 4, see [HW].

## 5 Algorithms For Computing Representations
The actual representations of numbers as sums of squares have various theoretical and practical applications. A repertoire of different algorithms have been developed for these tasks.

**6.1 Sums of Two Squares.**
Apparently there is no efficient algorithm for computing the representations of arbitrary integers as sums of two squares or the number of representations $r_2(N) = \sum_{d|N} \chi(d) = 4(d_1(N) - d_3(N))$, where $d_t(N) = \#\{ d \mid N : d \equiv t \bmod 4 \}$, as sums of two squares. However, there is an efficient algorithm for computing the representations of prime numbers.

***Theorem 7.*** (Hermite-Smith algorithm) The determination of the representation $p = x^2 + y^2$ has probabilistic logarithm time complexity.

A simplified version of the classical Hermite-Smith algorithm was derived in [BR]. There are other similar algorithms such as the Cornacchia algorithm, for computing the more general representations of primes $p = x^2 + dy^2$, with $|d| \geq 1$, see [CE, p. 35]. The determinations of the solutions for composite numbers $N = ax^2 + bxy + cy^2$ are significantly more difficult, see [HD].

**6.2 Sums of Three and Four Squares.**
Various schemes for computing the representations of integers as sums of three and four squares are provided and analyzed in [RS].

***Theorem 8.*** ([RS]) There is a randomized algorithm for expressing a number $N$ as a sum of four squares which requires an expected number of $O(\log^2(N)\log\log(N))$ operations with integers smaller than $N$ for all $N \geq N_0$.

*Sketch of the proof:* The number of solutions of the equation $N = x^2 + y^2 + p$, where $p = 4a + 1$ is a prime, is $cN/(\log(N)\log\log(N))$, $c > 0$ constant, this is Linnik Theorem. Thus, the expected number of random pairs $x, y \leq N^{1/2}$ before a solution is found is $O(\log(N)\log\log(N))$. Lastly, obtain $N = x^2 + y^2 + z^2 + w^2$ after replacing the representation of the prime $p = 4a + 1 = z^2 + w^2$. ∎

***Corollary 9.*** There is a probabilistic algorithm for determining the representation $N = x_1^2 + x_2^2 + \cdots + x_{2s}^2$ of an integer $N$ as a sum of $2s$ squares, $s \geq 1$. The algorithm has probabilistic logarithmic time complexity $O(\log^c(N))$, $c > 0$ constant, for all $N \geq N_0$.

*Proof*: Take $4sN = 2N + \cdots + 2N$, $s$ copies of $2N$, and apply Theorem 8 to each term $2N$. ∎

Apparently, there is no algorithm for computing the number of representations $r_4(N) = 8(d_1(N) + d_2(N) + d_3(N)) = 8\sum_{d|N, d \equiv 1,2,3 \bmod 4} d$ as sums of four squares.

# 6. REFERENCES

[BR] Brillhart, John. Note on representing a prime as a sum of two squares. Math. Comp. 26 (1972), 1011-1013.

[DN] Dickson, Leonard Eugene. History of the theory of numbers. Vol. I: Divisibility and primality. Chelsea Publishing Co., New York 1966.

[CE] Cohen, Henri. A course in computational algebraic number theory. Graduate Texts in Mathematics, 138. Springer-Verlag, Berlin, 1993.

[ES] Erdos, Paul; Nathanson, Melvyn B. Lagrange's theorem and thin subsequences of squares. Contributions to probability, pp. 3-9, Academic Press, New York-London, 1981.

[GR] Grosswald, Emil. Representations of integers as sums of squares. Springer-Verlag, New York, 1985. ISBN: 0-387-96126-7

[HD] Hardy, Kenneth; Muskat, Joseph B.; Williams, Kenneth S. Solving $n=au^2+buv+cv^2$n=au2+buv+cv2 using the Euclidean algorithm. Utilitas Math. 38 (1990), 225–236.

[HN] Glyn Harman, Small Fractional Parts Of Additive Forms In Prime Variables, Q J Math (1995) 46(3): 321-332

[HW] G. H. Hardy and E. M. Wright, An Introduction to the Theory of Numbers, 5th ed., Oxford University Press, Oxford, 1979.

[MW] Moreno, Carlos J.; Wagstaff, Samuel S., Sums of squares of integers. Discrete Mathematics and its Applications. Chapman & Hall, CRC, Boca Raton, FL, 2006. ISBN: 978-1-58488-456-9

[NT] Nathanson, Melvyn B. Elementary methods in number theory. Graduate Texts in Mathematics, 195. Springer-Verlag, New York, 2000. ISBN: 0-387-98912-9.

[RS] M. O. Rabin, J. O. Shallit, *Randomized Algorithms in Number Theory*, Communications on Pure and Applied Mathematics 39 (1986), no. S1, pp. S239–S256.

[ZR] Zöllner, Joachim. Über eine Vermutung von Choi, Erdos und Nathanson. Acta Arith. 45 (1985), no. 3, 211-213.